\documentclass{amsart}

\usepackage{amssymb}
\usepackage{amsmath}

\renewcommand{\arraystretch}{1.4}

\newtheorem{theorem}{Theorem}[section]
\newtheorem{lemma}[theorem]{Lemma}
\newtheorem{proposition}[theorem]{Proposition}
\newtheorem{corollary}[theorem]{Corollary}

\newtheorem{_definition}[theorem]{Definition}
\newenvironment{definition}{\begin{_definition}\rm}{\end{_definition}}

\newtheorem{_remark}[theorem]{\it Remark}
\newenvironment{remark}{\begin{_remark}\rm}{\end{_remark}}

\newtheorem{_claim}[theorem]{Claim}
\newenvironment{claim}{\begin{_claim}\rm}{\end{_claim}}

\numberwithin{equation}{section}
\numberwithin{table}{section}
\numberwithin{figure}{section}

\renewcommand{\P}{\mathord{\mathbb  P}}
\newcommand{\Z}{\mathord{\mathbb Z}}
\newcommand{\R}{\mathord{\mathbb R}}
\newcommand{\F}{\mathord{\mathbb  F}}

\newcommand{\Q}{\mathord{\mathbb Q}}

\newcommand{\CCC}{\mathord{\mathcal C}}

\newcommand{\EEE}{\mathord{\mathcal E}}

\newcommand{\LLL}{\mathord{\mathcal L}}
\newcommand{\OOO}{\mathord{\mathcal O}}

\newcommand{\maprightsp}[1]{\smash{\mathop{\; \longrightarrow \; }\limits\sp{#1}}}

\newcommand{\inj}{\mathbin{\hookrightarrow}}

\newcommand{\isom}{\smash{\mathop{\;\to\;}\limits\sp{\sim}}}

\newcommand{\set}[2]{\{\; {#1} \; \mid \; {#2} \;  \}}
\newcommand{\shortset}[2]{\{ {#1}  \mid  {#2}  \}}

\newcommand{\map}[3]{ #1 \, : \, #2 \, \to \, #3}

\newcommand{\sm}{\setminus}
\newcommand{\st}{\subset}

\newcommand{\ol}[1]{\overline{#1}}

\newcommand{\sprime}{\sp\prime}
\newcommand{\sptimes}{\sp\times}
\newcommand{\spprime}{\sp{\prime\prime}}

\newcommand{\dual}{\sp{\vee}}

\newcommand{\sperp}{\sp\perp}

\newcommand{\inv}{\sp{-1}}

\newcommand{\rmand}{\textrm{and}}

\newcommand{\quand}{\quad\rmand\quad}

\newcommand{\Pt}{\P^2}

\newcommand{\parder}{\partial}

\newcommand{\Spec}{\operatorname{\mathord{\rm{Spec}}}\nolimits}
\newcommand{\Proj}{\operatorname{\mathord{\rm{Proj}}}\nolimits}
\newcommand{\NS}{\mathord{\rm{NS}}}
\newcommand{\ch}{\operatorname{\mathord{\rm{char}}}\nolimits}
\newcommand{\disc}{\operatorname{\mathord{\rm{disc}}}\nolimits}
\newcommand{\Roots}{\mathord{\rm{Roots}}}
\newcommand{\Nef}{\mathord{\rm{Nef}}}

\newcommand{\Hom}{\mathord{\rm{Hom}}}
\newcommand{\Sing}{\mathord{\rm{Sing}}}
\newcommand{\rank}{\operatorname{\mathord{\rm{rank}}}\nolimits}

\newcommand{\pol}{\mathord{H}}

\newcommand{\Pp}{\mathord{\mathcal{P}}}
\newcommand{\Qp}{\mathord{\mathcal{Q}}}
\newcommand{\SP}{\mathord{\mathcal{L}}}

\newcommand{\lat}{\Lambda}
\newcommand{\dgr}{\Delta}

\newcommand{\angl}[1]{\langle #1\rangle}

\newcommand{\rs}{_{r, s}}
\newcommand{\ab}{\alpha\beta}

\newcommand{\ve}{\varepsilon}

\begin{document}
\title[Supersingular $K3$ surfaces in characteristic $2$]{
On Kummer type construction of supersingular $K3$ surfaces in characteristic $2$}

\author{Ichiro Shimada}
\address{
Division of Mathematics,
Graduate School of Science,
Hokkaido University,
SAPPORO 060-0810, JAPAN
}
\email{shimada@math.sci.hokudai.ac.jp
}

\author{De-Qi Zhang}
\address{
Department of Mathematics,
National University of Singapore,
Lower Kent Ridge Road, 
Singapore 119260
}
\email{matzdq@math.nus.edu.sg
}

\begin{abstract}
We show that every supersingular $K3$ surface in characteristic $2$
with Artin invariant $\le 2$ is obtained by the Kummer type construction of Schr\"oer.
\end{abstract}


\maketitle

\section{Introduction}\label{sec:intro}
We work over an algebraically closed field $k$.
A $K3$ surface $X$ is called \emph{supersingular} (in the sense of Shioda)
if the rank of the N\'eron-Severi lattice $\NS (X)$ of $X$ attains the possible maximum $22$.
Supersingular $K3$ surfaces exist only when $\ch k$ is positive.
The Artin invariant $\sigma (X)$ of a supersingular $K3$ surface $X$
 is defined in~\cite{Artin1975} by
$$
\disc \NS (X)=-p^{2\sigma (X)},
$$
where $p=\ch k>0$.
It is known that $\sigma (X)$ is a positive integer $\le 10$.

Let $A$ be an abelian surface,
and let $\iota :A\to A$ be the involution $x\mapsto -x$.
If $\ch k\ne 2$,
then the minimal resolution of the quotient surface $A/\langle \iota\rangle$
is a $K3$ surface,
which is called the \emph{Kummer surface associated with $A$}.

An abelian surface $A$ in positive characteristic is called \emph{supersingular}
if $A$ is isogenous to a product of supersingular elliptic curves.
Ogus~\cite{MR563467, Ogus}
 proved that,
if $\ch k>2$,
the supersingular $K3$ surfaces with Artin invariant $\le 2$
are exactly the Kummer surfaces associated with supersingular abelian surfaces.
(See also Shioda~\cite{Shioda1}.)
On the other hand,
Shioda~\cite{Shioda2} and Katsura~\cite{Katsura} observed  that,
if $\ch k=2$, 
then the minimal resolution of   the quotient
of a supersingular abelian surface  by the involution $x\mapsto -x$ 
is  a rational surface.

In~\cite{Sch}, Schr\"oer presented a  Kummer type construction of 
supersingular $K3$ surfaces in characteristic $2$.
We assume that $\ch k=2$ in this paragraph.
Let $C\times C$ be the self-product of the rational curve $C$ 
with one ordinary cusp.
We put
$$
\begin{array}{cccl}
C &=& \Spec k [u^2, u^3]\cup \Spec k[u^{-1}] &\textrm{for the first factor, and}\cr
C &=& \Spec k [v^2, v^3]\cup \Spec k[v^{-1}] &\textrm{for the second factor.}
\end{array}
$$
Let $r$ and $s$ be constants in $k$ such that $(r, s)\ne (0,0)$.
Then the derivation
\begin{equation}\label{eq:delta}
 (u^{-2} +r) \frac{\parder}{\parder u}+(v^{-2} +s)\frac{\parder}{\parder v}
\end{equation}
defines a global vector  field $\delta$ on $C\times C$ satisfying $\delta^{[2]}=0$.
Hence $\delta$ corresponds to an action of 
the infinitesimal  group scheme $\alpha_2$ on $C\times C$.
Let $X_{r,s}$ be the minimal resolution of the quotient surface 
$(C\times C)/\alpha_2$.
\begin{theorem}[\cite{Sch}]
The surface $X_{r, s}$ is a supersingular $K3$ surface with Artin invariant
$$
\sigma (X_{r,s})=\begin{cases}
1 & \textrm{if $r=0$ or $s=0$ or $r^3=s^3$,}\\
2 & \textrm{otherwise}.
\end{cases}
$$
\end{theorem}
The purpose of this paper is to prove the following:
\begin{theorem}\label{thm:main}
Let $X\sprime$ be a supersingular $K3$ surface 
in characteristic $2$ with Artin invariant $\le 2$.
Then there exist  constants $r,s\in k$ with $(r, s)\ne (0,0)$ such that
$X\sprime $ is isomorphic to  Schr\"oer's Kummer surface $X_{r, s}$.
\end{theorem}
Even though the moduli curve of  marked supersingular $K3$ surfaces
with Artin invariant $\le 2$ is constructed~(\cite{Ogus, RS}), 
it is not separated.
Hence the existence of the complete  family  of Schr\"oer's Kummer surfaces
of dimension $1$
does not imply Theorem~\ref{thm:main} immediately.
\par
\medskip
The main ingredient of the proof is the following 
structure theorem for N\'eron-Severi lattices of supersingular $K3$ surfaces
due to  Rudakov and Shafarevich~\cite{RS}:
\begin{theorem}\label{thm:RS}
Let $X$ and $X\sprime$ be supersingular $K3$ surfaces 
defined over the  same algebraically closed field.
If $\sigma (X)=\sigma (X\sprime)$, then
the lattices $\NS (X)$ and $\NS (X\sprime)$ are isomorphic.
\end{theorem}
Indeed, the N\'eron-Severi lattice $\NS (X)$ of a supersingular $K3$ surface $X$
in characteristic $p$ is $p$-elementary~(\cite[Theorem in Section 8]{RS}, see also~\cite{Artin1975}).
If $p=2$, then $\NS (X)$ is of type I~(\cite[Proposition in Section 5]{RS}).
Hence the classification theorem of even 
hyperbolic $p$-elementary lattices~(\cite[Theorem in Section 1]{RS})
implies Theorem~\ref{thm:RS}.
\par
\medskip
The outline of the proof of Theorem~\ref{thm:main} is as follows.
First note that,
by~\cite[Proposition 6.2]{Sch},
if $\sigma (X\rs)=2$, then   Schr\"oer's Kummer surface 
$X\rs$ is birational to a purely inseparable
double cover 
$Y\rs$ of $\Pt$ defined by
$$
w^2=x(y^4+s^2y^2)+y(x^4+r^2x^2),
$$
which has rational double points of type $4D_4+5A_1$.
Let us assume, for simplicity,  that 
the given supersingular $K3$ surface $X\sprime$ is of Artin invariant $2$. 
We choose one of  Schr\"oer's Kummer surfaces $X$ with Artin invariant $2$
(for example, we put $X:=X_{1,s}$ with $s\notin \F_4$).
Using the isomorphism between $\NS (X)$ and $\NS (X\sprime)$,
we can show that $X\sprime$ is also birational to a double cover $Y\sprime$
of $\Pt$ with rational double points of  type $4D_4+5A_1$.
By means of the 
 notion of \emph{half-lines} and \emph{splitting lines},
we can show that the covering morphism $Y\sprime\to \Pt$ is purely inseparable,
and  then we can determine the defining equation of $Y\sprime$.
It turns out that the defining equation of $Y\sprime$ is equal to that of $Y_{t, 1}$
for some non-zero constant $t\in k$.
Therefore $X\sprime$ is isomorphic to  Schr\"oer's Kummer surface $X_{t, 1}$.
\par
\medskip
A surface birational to a purely inseparable  cover of $\Pt$
is called a \emph{Zariski surface},
and its basic properties have been studied in~\cite{BL}.
In~\cite{Shimada1}~and~\cite{Shimada2},
we showed that every supersingular $K3$ surface in characteristic $2$
is birational to a purely inseparable double cover of $\Pt$
with $21$ ordinary nodes,
and studied the N\'eron-Severi lattice of such a surface.
Using the results obtained in~\cite{Shimada2},
we have determined in~\cite{Shimada3}
the moduli curve of polarized supersingular $K3$ surfaces
with Artin invariant $\le 2$ and with $21$ ordinary nodes.
In~\cite{Sch}, Schr\"oer showed that, as $r$ and $s$ varies, 
his Kummer surfaces $X\rs$ form a smooth family over the projective line
$\Proj k[\sqrt{r}, \sqrt{s}]$.
It would be an interesting problem to investigate the relation between
the moduli curve in~\cite{Shimada3} and Schr\"oer's projective line.
\par
\medskip
On the other hand, in~\cite{ShimadaZhang}, 
we investigated 
supersingular $K3$ surfaces with $10$ ordinary cusps.
Such  supersingular $K3$ surfaces exist only in characteristic $3$.
An example is obtained as a purely inseparable \emph{triple} cover of $\P^1\times\P^1$.
The proof in the present article of the fact that  $Y\sprime\to \Pt$ is purely inseparable
uses an argument developed in~\cite{ShimadaZhang}.
\par
\medskip
The plan of this paper is as follows.
In~\S\ref{sec:lattice},
we collect from the lattice theory
some definitions and facts
that will be used in this paper.
The very elementary Lemmas~\ref{lem:A}~and~\ref{lem:D}
play an important role in the proof 
of the fact that $Y\sprime$ is purely inseparable over $\Pt$.
In~\S\ref{sec:polK3},
we review some properties of the N\'eron-Severi lattice of a $K3$ surface.
We then introduce the notion of half-lines and splitting lines
for a polarized $K3$ surface of degree $2$
in~\S\ref{sec:pol2}.
After investigating 
the purely inseparable double cover $Y\rs\to\Pt$
birational to 
 Schr\"oer's Kummer  surface $X\rs$,
we prove Theorem~\ref{thm:main} in~\S\ref{sec:proof}.
\section{Preliminaries on lattices}\label{sec:lattice}
A free $\Z$-module $\lat$ of finite rank with a non-degenerate
symmetric bilinear form
\begin{equation}\label{eq:billat}
\lat\times \lat \to \Z
\end{equation}
denoted by $(u, v)\mapsto uv$ is called a \emph{lattice}.
Let $\lat$ be a lattice.
The \emph{dual lattice} $\lat\dual$ of $\lat$ is the $\Z$-module 
$\Hom (\lat, \Z)$.
Then $\lat$ is naturally embedded into $\lat\dual$ as a submodule of finite index.
The \emph{discriminant group} 
of $\lat$ is, by definition,
the finite abelian group $\lat\dual/\lat$.
There exists a unique symmetric bilinear form
\begin{equation}\label{eq:billatdual}
\lat\dual\times \lat\dual \to \Q
\end{equation}
that extends~\eqref{eq:billat}.
An \emph{overlattice} of $\lat$ is a submodule $N$ of $\lat\dual$
containing $\lat$
such that the bilinear form~\eqref{eq:billatdual}
takes values in $\Z$ on $N\times N$.
If $\lat$ is a sublattice of a lattice $\lat\sprime$ with finite index,
then $\lat\sprime$ is 
embedded into $\lat\dual$  in a natural way,
and hence is  regarded as an overlattice of $\lat$.
\par
We say that $\lat$ is \emph{even} if $u^2\in 2\Z$ holds for every $u\in \lat$.
The signature $(s_+, s_-)$ of a lattice $\lat$ is the numbers of positive and negative eigenvalues of
the intersection matrix of $\lat$.
We say 
that $\lat$ is \emph{negative-definite}
if $s_+=0$, 
and 
that $\lat$ is \emph{hyperbolic}
if $s_+=1$.
By abuse of language, 
a positive definite lattice of rank $1$ is also called hyperbolic.
\par
Let $\lat$ be an even negative-definite lattice.
A vector $r\in \lat$ is called a \emph{root}
if $r^2=-2$.
We denote by $\Roots (\lat)$ the set of roots in $\lat$.
We define an equivalence relation $\sim$ on $\Roots (\lat)$
by the following:
$r\sim r\sprime$ if there exists a sequence $r_0=r, r_1, \dots, r_{m-1},  r_m=r\sprime$
of roots in $\lat$ such that $r_i r_{i+1}\ne 0$ for $i=0, \dots, m-1$.
Let $R_1, \dots, R_k$ be the equivalence classes of $\sim$.
We call the decomposition
$$
\Roots (\lat)=R_1\sqcup\dots\sqcup R_k
$$
the \emph{irreducible decomposition} of $\Roots (\lat)$.
Suppose that we are given a linear form
$$
\map{\alpha}{ \lat}{\R}
$$
such that $\alpha (r)\ne 0$ for any $r\in \Roots (\lat)$.
We put
\begin{equation}\label{eq:Riplus}
R_i\sp+ := \set{r\in R_i}{\alpha (r)>0}.
\end{equation}
A root $r\in R_i\sp{+}$ is called \emph{decomposable} if there exist
$r_1, r_2\in  R_i\sp{+}$ such that $r=r_1+r_2$,
and $r$ is called \emph{indecomposable} if it is not decomposable.
For the proof of the following results,
see~\cite{Ebeling} or~\cite{Bourbaki}, for example.
\begin{proposition}\label{prop:nonnegative}
Let $r$ be an element of $R_i\sp{+}$ such that $\alpha (r)>0$.
Then $r$ can be  written in a unique way as a linear combination of 
indecomposable elements of $R_i\sp{+}$.
Moreover the coefficients are all non-negative integers.
\end{proposition}
\begin{proposition}\label{prop:latroots}
Let $\lat_i$ be the sublattice of $\lat$ generated by the roots in $R_i$.
Then $\lat_1, \dots, \lat_k$ form an orthogonal direct sum in $\lat$.
The indecomposable elements of $R_i\sp+$ form a basis of the lattice $\lat_i$,
and the intersection matrix of $\lat_i$
with respect to this basis is a  Cartan matrix
of type $ADE$ multiplied by $-1$. 
\end{proposition}
The indecomposable elements of $R_i\sp+$ have the following characterization:
\begin{corollary}\label{cor:indecomp}
Let $\ve_1, \dots, \ve_d$ be elements of $R_i\sp+$ such that
every element of $R_i\sp+$ is written as a linear combination of
$\ve_1, \dots, \ve_d$ 
with non-negative integer coefficients in a unique way.
Then $\{\ve_1, \dots, \ve_d\}$ is equal to the set of 
indecomposable elements of $R_i\sp+$.
\end{corollary}
\begin{proof}
Suppose that $\ve_i$ is decomposable.
There exist $r_1, r_2\in R_i\sp+$ such that $\ve_i=r_1+r_2$.
Since each of $r_1$ and $r_2$ is written 
as a linear combination of
$\ve_1, \dots, \ve_d$ 
with non-negative integer coefficients,
we obtain a contradiction to the uniqueness of the way to write
$\ve_i$ as a linear combination of
$\ve_1, \dots, \ve_d$ 
with non-negative integer coefficients.
Therefore each of $\ve_1, \dots, \ve_d$ is indecomposable.
\par
Suppose that $r\in R_i\sp+$ is indecomposable.
We can write $r$ as a linear combination of
$\ve_1, \dots, \ve_d$ 
with non-negative integer coefficients.
Since each $\ve_i$ is indecomposable,
the uniqueness of the way to write
$r$ as a 
linear combination of
indecomposable elements of $R_i\sp+$
with non-negative integer coefficients
implies that $r$ is equal to one of $\ve_1, \dots, \ve_d$.
\end{proof}
Let $\tau_i$ be the $ADE$-type of the Cartan matrix of the intersection matrix of $\lat_i$ given in Proposition~\ref{prop:latroots}.
We define the \emph{root type of $\lat$} to be the formal sum
$\tau_1+\cdots+ \tau_k$.
\par
\medskip
We say that $\lat$ is a \emph{root lattice}
if $\lat$ is generated by $\Roots(\lat)$.
For later use, we present properties of root lattices of  type
$A_1$ and $D_4$.
\par
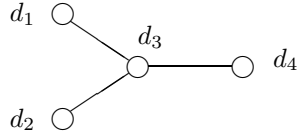
\begin{figure}[h]
\begin{center}
\setlength{\unitlength}{1mm}
\begin{picture}(50, 16)(0,0)
\put(-2, 14){$d_1$}\put(5,15){\circle{2.8}}
\put(-2, 0){$d_2$}\put(5,1){\circle{2.8}}
\put(15, 11){$d_3$}\put(15,8){\circle{2.8}}
\put(33, 8){$d_4$}\put(29,8){\circle{2.8}}
\put(6.1, 14){\line(3, -2){7.7}}
\put(6.1, 2.0){\line(3, 2){7.7}}
\put(16.5, 8){\line(1,0){11}}
\end{picture}
\end{center}
\caption{The Dynkin diagram of type $D_4$}\label{fig:dynkinD4}
\end{figure}

Let $\lat$ be the root lattice of type $A_1$,
and let $a\in \lat$ be a root,
which generates $\lat$.
We put $a\dual:=-a/2$,
which generates $\lat\dual$.
Then the discriminant group of $\lat$ is isomorphic to $\Z/2\Z$.
The proof of the following is elementary:
\begin{lemma}\label{lem:A}
Let $v\in \lat\dual$ be a vector such that $va\ge 0$.
If $v\equiv 0\bmod \lat$,
then we have $v^2=0$ or $v^2\le -2$,
and $v^2=0$ holds if and only if $v=0$.
If $v\equiv a\dual\bmod \lat$,
then we have $v^2=-1/2$ or $v^2\le -9/2$,
and $v^2=-1/2$ holds if and only if $v=a\dual$.
\end{lemma}
Let $\lat$ be the root lattice of type $D_4$
generated by the roots $d_1, \dots, d_4$
whose intersection numbers  are given by the Dynkin diagram  in Figure~\ref{fig:dynkinD4}.
Let $d_1\dual, \dots, d_4\dual$ be the basis of $\lat\dual$ 
dual to $d_1, \dots, d_4$.
We have
\begin{equation}\label{eq:dualD4}
\renewcommand{\arraystretch}{1.4}
[d_1\dual, d_2\dual, d_3\dual, d_4\dual]=
[d_1, d_2, d_3, d_4]
\left[
\begin{array}{cccc}
-1 &-1/2 & -1 &-1/2\\
-1/2 & -1 & -1 & -1/2 \\
-1 & -1 & -2 & -1\\
-1/2 & -1/2 & -1 & -1
\end{array}
\right].
\end{equation}
The discriminant group of $\lat$ is isomorphic to $(\Z/2\Z)\oplus (\Z/2\Z)$,
and is generated by $d_1\dual \bmod \lat$ and $d_4\dual\bmod \lat$.
\begin{lemma}\label{lem:D}
Let $v\in\lat\dual$ be a vector such that
$vd_i\ge 0$ holds for $i=1, \dots, 4$.
If $v\equiv 0\bmod \lat$,
then we have $v^2=0$ or $v^2\le -2$,
and $v^2=0$ holds if and only if $v=0$.
If $v\equiv d_1\dual \bmod\lat$,
then we have $v^2=-1$ or $v^2\le -3$,
and $v^2=-1$ holds if and only if $v=d_1\dual$.
\end{lemma}
\begin{proof}
The first assertion is obvious.
Suppose that  $v\equiv d_1\dual\bmod \lat$.
Then we can put
$$
v=d_1\dual + x_1 d_1 + x_2 d_2 +x_3 d_3 + x_4 d_4,
$$
where $x_1, \dots, x_4\in \Z$.
From the condition $vd_i\ge 0$ for $i=1, \dots, 4$,
we obtain the following inequalities:
\begin{equation}\label{eq:ineq}
1-2x_1+x_3\ge 0, 
\;
-2x_2 + x_3\ge 0, 
\;
x_1 + x_2 -2 x_3 +x_4\ge 0, 
\;
x_3-2x_4\ge 0.
\end{equation}
Using~\eqref{eq:dualD4},
we calculate 
$$
\begin{array}{ccl}
v^2 &=& -1 -2 (x_1^2 +x_2^2 +x_3^2 +x_4^2 -x_1 -x_1 x_3 - x_2 x_3 - x_3 x_4)\\
&=& -1 -\{(1-2x_1+x_3)^2+(-2x_2 + x_3)^2+(x_3-2x_4)^2+(x_3-1)^2-2\}/2.
\end{array}
$$
Therefore $v^2$ is a negative odd integer,
and $v^2=-1$ holds if and only if
two of the four integers $1-2x_1+x_3$, $-2x_2 + x_3$, $x_3-2x_4$, $x_3-1$ are $\pm 1$
and the other two are $0$.
Combining this with the inequalities~\eqref{eq:ineq},
we see that 
$v^2=-1$ holds if and only if $x_1=x_2=x_3=x_4=0$.
\end{proof}
\section{The N\'eron-Severi lattice of a $K3$ surface}\label{sec:polK3}
In this section,
we work over an algebraically closed field 
of arbitrary characteristic.
Let $X$ be an (algebraic) $K3$ surface,
and let $\NS (X)$ be the N\'eron-Severi lattice of $X$,
which is an even hyperbolic lattice.
For a divisor $D$ on $X$,
we denote by $[D]\in \NS(X)$ the class of $D$.
%
%
\subsection{The nef-cone}\label{subsec:nef}
We put 
$$
\Nef (X):=\set{v\in \NS (X)\otimes \R}{\textrm{%
$v[D]\ge 0$ for any effective divisor  $D$ on $X$}}.
$$
Let $A$ be an ample divisor on $X$,
and let $\CCC^+ (X)$ be the connected component of
$$
\set{v\in \NS (X)\otimes \R}{v^2>0}
$$
that contains $[A]$.
For a vector $v\in \NS(X)$, we put
$$
\angl{v}_{\R}\sperp:=\set{w\in \NS (X)\otimes \R}{vw=0}.
$$
Then the family of hyperplanes
$\shortset{\angl{r}_{\R}\sperp}{r^2=-2}$ 
of $\NS (X)\otimes \R$ is locally finite in $\CCC^+ (X)$.
It is well-known and easy to prove that
$[A]\notin \angl{r}_{\R}\sperp$ for any vector $r$ with $r^2=-2$, 
and that
$\Nef (X)$ is equal to the closure in $\NS (X)\otimes \R$
of the connected component of 
$$
\CCC^+ (X)\; \setminus\; \bigcup \;\angl{r}_{\R}\sperp
$$
that contains $[A]$.
By the argument of Proposition 3 in~\cite[Section 3]{RS}, we obtain the following:
\begin{proposition}\label{prop:nef}
Let $X$ and $X\sprime$ be two algebraic $K3$ surfaces such that $\NS (X)$ and $\NS(X\sprime)$ are
isomorphic. Then there exists an isomorphism
$\phi: \NS (X)\isom \NS(X\sprime)$ such that $\phi\otimes\R$
maps $\Nef (X)$ to $\Nef (X\sprime)$.
\end{proposition}
\subsection{Polarizations}\label{subsec:pol}
\begin{proposition}\label{prop:nikulin}
Let $H$ be a divisor on an algebraic $K3$ surface 
 $X$ such that  $[H]\in \Nef (X)$ and $H^2>0$.
Then the following conditions are equivalent to each other:
\begin{itemize}
\item[(i)] The complete linear system $|H|$  has no fixed components.
\item[(ii)]  There exist no vectors $e\in \NS (X)$
such that $e[H]=1$ and $e^2=0$.
\end{itemize}
\end{proposition}
\begin{proof}
The implication (i)$\Longrightarrow$(ii) follows from the argument in 
the proof of (4)$\Longrightarrow$(1) in~\cite[Proposition~1.7]{Urabe}.
The other implication (ii)$\Longrightarrow$(i) follows from~\cite[Proposition 0.1]{Nikulin}.
\end{proof}
\begin{definition}
A  \emph{polarization} of an algebraic $K3$ surface $X$ is 
 a divisor $H$ on $X$ satisfying 
$[H]\in \Nef (X)$, $H^2>0$, and the conditions (i) and (ii) 
in Proposition~\ref{prop:nikulin}.
The positive integer $H^2$ is called the \emph{degree}
of the polarization $H$.
By~\cite[Proposition~0.1]{Nikulin}, if $H$ is a polarization of degree $d$, then 
$|H|$ is base-point free by  Saint-Donat~\cite[Corollary 3.2]{MR0364263}, and  we have
$\dim |H|=1+d/2$.

A pair $(X, H)$ of a $K3$ surface $X$ and a polarization $H$ of 
$X$ is called a \emph{polarized $K3$ surface}.
\end{definition}
Combining Propositions~\ref{prop:nef}~and~\ref{prop:nikulin}, we obtain the following:
\begin{corollary}\label{cor:pol}
Let $X$ and $X\sprime$ be two $K3$ surfaces such that $\NS (X)$ and $\NS(X\sprime)$ are isomorphic,
and
let $H$ be a polarization of $X$.
If $\phi: \NS (X)\isom \NS(X\sprime)$ is an isomorphism 
 such that $\phi\otimes\R$
maps $\Nef (X)$ to $\Nef (X\sprime)$,
then $\phi([H])$ is the class of a polarization $H\sprime$ of $X\sprime$.
\end{corollary}
A curve $C$ on $X$ is called a \emph{$(-2)$-curve} on $X$ if 
it satisfies the following conditions that are equivalent to each other:
\begin{itemize}
\item[(i)] $C$ is a smooth rational curve,
\item[(ii)] $C$ is reduced irreducible with negative self-intersection,
\item[(iii)] $C$ is irreducible and $C^2=-2$.
\end{itemize}
Let $(X, H)$ be a polarized $K3$ surface.
Then the complete linear system $|H|$ defines a morphism $\Phi_{|H|}$ from 
$X$ to a projective space $\P^N$ $(N=1+H^2/2)$
that is generically finite over the image.
We denote by
\begin{equation}\label{eq:stein}
X\;\;\maprightsp{\rho}\;\; Y\;\; \maprightsp{\pi}\;\; \P^N
\end{equation}
the Stein factorization of  $\Phi_{|H|}$;
that is, $\rho$ is birational, $Y$ is normal, and $\pi$ is finite. 
The normal $K3$ surface $Y$ has only rational
double points as its singularities,
and hence $\rho$ is a contraction of  
an $ADE$-configuration of $(-2)$-curves. 
(See~\cite{Artin1962, Artin1966}.)
Let $\EEE$ be the set of $(-2)$-curves that are contracted by $\rho$.
The classes $[E]$ of $E\in\EEE$ are determined by the following
procedure.
Let $[H]\sperp$ be  the orthogonal complement 
of $[H]$ in $ \NS (X)$.
Since $\NS (X)$ is even hyperbolic and $[H]^2$ is positive,
$[H]\sperp$ is even and negative-definite.
We can therefore consider the set $\Roots ([H]\sperp)$
of roots in $[H]\sperp$.
\begin{lemma}\label{lem:DE}
Let $r$ be an element of $\Roots ([H]\sperp)$.
Then there exists a unique effective divisor $E$
such that $r=[E]$ or $r=-[E]$ holds.
Moreover, the integral  component of $E$ is  a $(-2)$-curve.
\end{lemma}
\begin{proof}
By the Riemann-Roch theorem and the Serre duality,
we see that either $r$ or $-r$ is the class of an effective divisor.
Replacing $r$ with $-r$, if necessary,
we can assume that $r$ is the class of an effective divisor $E$.
Let $E=F+M$ be the decomposition of $E$ into the sum
of the fixed part $F$ and the movable part $M$.
Since $[H]\in \Nef (X)$,
we have $HF\ge 0$ and $HM\ge 0$.
Because $HE=0$, we have $HM=0$.
Since $[H]\sperp$ is negative-definite
and $M^2\ge 0$, we obtain 
 $M=0$.
Therefore $E$ is unique and every irreducible component
of $E$ has negative self-intersection number.
Thus the reduced part of every irreducible component
of $E$  is  a $(-2)$-curve.
\end{proof}
Let $\Roots ([H]\sperp)=R_1\sqcup\dots\sqcup R_k$
be the irreducible decomposition of $\Roots ([H]\sperp)$
defined in \S\ref{sec:lattice}.
We choose an interior point $a$ of $\Nef (X)$
(for example, the class of an ample divisor on $X$),
and let
$\alpha: \NS (X)\to \R$ be the linear form given by 
$\alpha (x):=ax$.
By Lemma~\ref{lem:DE},
we see that $\alpha (r)\ne 0$ for any $r\in \Roots([H]\sperp)$.
We thus can define $R_i\sp+$ by~\eqref{eq:Riplus},
and consider the indecomposable roots of $R_i\sp+$.
Note that $R_i\sp+ \subset R_i$ does not depend on the choice of 
the interior point $a$ of $\Nef (X)$.
\begin{proposition}
Let $\Sing (Y)$ be the set of singular points of $Y$.
There exists a bijection 
from the set $\{R_1, \dots, R_k\}$ to
$\Sing (Y)$
with the following property.
Let $P_i\in \Sing (Y)$ be the point corresponding to $R_i$.
Then the classes of $(-2)$-curves \ contracted by 
$\rho$ to  $P_i$ 
are exactly the indecomposable roots of $R_i\sp+$.
\end{proposition}
\begin{proof}
Let $r$ be an element of $R_i\sp+$.
By Lemma~\ref{lem:DE} and $\alpha (r)>0$,
$r$ is the class of a unique effective divisor  of the form
$a_1 E_1+\cdots+a_l E_l$,
where $E_1, \dots, E_l$ are $(-2)$-curves and $a_1, \dots, a_l$ are positive integers.
Since $[H]\in \Nef(X)$ and $r\in [H]\sperp$,
we have $[E_\nu]\in [H]\sperp$ for $\nu=1, \dots, l$.
In particular, we have  $E_\nu\in \EEE$  for $\nu=1, \dots, l$.
Let $\lat_j$ be the sublattice of $[H]\sperp$ generated by the roots in $R_j$ for $j=1, \dots, k$.
Since $\lat_1, \dots, \lat_k$ form a direct sum in $\NS (X)$,
the uniqueness of the effective divisor representing  $r\in \lat_i$
implies that $[E_1], \dots, [E_l]$ are all in $R_i$.
Since $\alpha ([E_\nu])>0$, we have $[E_\nu]\in R_i\sp+$.
Thus we have shown that every element of $R_i\sp+$
is written as a linear combination of the classes of $(-2)$-curves in $R_i\sp+$
with non-negative integer coefficients in a unique way.
By Corollary~\ref{cor:indecomp}, we see that 
$r$ is the class of a $(-2)$-curve in $\EEE$ 
if and only if $r$ is  indecomposable in $R_i\sp+$.
\end{proof}
Let $(X\sprime, H\sprime)$ be another polarized $K3$ surface.
Let
\begin{equation*}
X\sprime\;\;\maprightsp{\rho\sprime}\;\; Y\sprime\;\; \maprightsp{\pi\sprime}\;\; \P^{N\sprime}
\end{equation*}
be the Stein factorization of the morphism $\Phi_{|H\sprime|}$
defined by $|H\sprime|$, and 
let $\EEE\sprime$ be the set of $(-2)$-curves contracted by $\rho\sprime$.
\begin{corollary}\label{cor:commutative}
Suppose that there exists an isomorphism 
$\phi: \NS(X)\isom\NS(X\sprime)$
such that $\phi\otimes \R$ maps $\Nef (X)$ to $\Nef (X\sprime)$,
and that $\phi ([H])$ is equal to $[H\sprime]$.
Then the $ADE$-type of $\Sing (Y)$ coincides with that of $\Sing (Y\sprime)$.
Moreover,
there exist bijections
$$
\phi_{\EEE}: \EEE\isom \EEE\sprime\quand
\phi_{\Sing}: \Sing (Y)\isom \Sing(Y\sprime)
$$
such that the following diagram is commutative;
\begin{equation}\label{eq:commutative}
\begin{array}{cccc}
\NS (X)&\maprightsp{\phi} & \NS(X\sprime)& \\
\big{\uparrow} && \big{\uparrow}& \\
\EEE &\maprightsp{\phi_{\EEE}} & \EEE\sprime&\\
\big{\downarrow} &&\big{\downarrow}&\\
\Sing (Y) & \maprightsp{\phi_{\Sing}} & \Sing(Y\sprime)&,\\
\end{array}
\end{equation}
where the up-arrows are given by
$E\mapsto [E]\in \NS(X)$ and $E\sprime\mapsto [E\sprime]\in \NS(X\sprime)$,
respectively,  and the down-arrows are given by
$E\mapsto \rho (E)\in\Sing (Y)$ and 
$E\sprime \mapsto \rho\sprime(E\sprime)\in\Sing (Y\sprime)$,
respectively.
\end{corollary}
\subsection{Polarizations with maximal rational double points}\label{subsec:polmax}
\begin{definition}
We say that a polarized $K3$ surface $(X, H)$ has \emph{maximal rational double points}
if the total Milnor number of $\Sing (Y)$ is equal to
$\rank \NS (X)-1$;
or equivalently,
the root lattice generated by $\Roots ([H]\sperp)$ is of finite index in $[H]\sperp$.
\end{definition}
Let $(X, H)$ be a polarized $K3$ surface with maximal rational double points.
Consider the Stein factorization~\eqref{eq:stein} of $\Phi_{|H|}$.
For $P\in \Sing (Y)$,
we denote by $\EEE_P\subset \EEE$ the set of $(-2)$-curves 
that are contracted to $P$ by $\rho$,
by $\lat_P\subset \NS (X)$
the sublattice generated by the classes $[E]$ of the curves $E\in \EEE_P$,
and by $\dgr_P$ the discriminant group $\lat_P\dual/\lat_P$ of $\lat_P$.
We also denote by $\lat_{\pol}\subset \NS (X)$ the sublattice of rank $1$ generated by $[H]$,
and by $\dgr_{\pol}$ the discriminant group $\lat_{\pol}\dual/\lat_{\pol}$ of $\lat_{\pol}$,
which is a cyclic group of order equal to  $H^2$.
We then put
$$
\lat:= \lat_{\pol} \oplus \bigoplus_{P\in \Sing (Y)}\lat_P
\quand
\dgr:=\lat\dual/\lat.
$$
We have natural decompositions
$$
\lat\dual= \lat_{\pol}\dual \oplus \bigoplus_{P\in \Sing (Y)}\lat_P\dual\quad\quand\quad\dgr= \dgr_{\pol} \oplus \bigoplus_{P\in \Sing (Y)}\dgr_P.
$$
By the assumption,
$\lat$ is of finite index in $\NS(X)$,
and hence $\NS(X)$ is an overlattice of $\lat$.
Let $v$ be a vector of $\NS (X)$.
Using the direct-sum decomposition of $\lat\dual$ and the natural embedding $\NS(X)\inj\lat\dual$,
we can  define the \emph{$H$-component} $v_H\in \lat\dual_{\pol}$  
and the \emph{$P$-components $v_P\in \lat\dual_{P}$} of $v$.
We denote by $\bar{v}\in \dgr$ the class of $v$ modulo $\lat$.
Then  the \emph{$H$-component} $\bar{v}_H\in \dgr_{\pol}$  
and the \emph{$P$-components $\bar{v}_P\in \dgr_{P}$} of $\bar{v}$ are also defined.
\section{Polarizations of degree $2$ in characteristic $2$}\label{sec:pol2}
From now on to the end of this paper,
we assume that the base field $k$ is of characteristic $2$.
\par
\medskip
Let $(X, H)$ be a polarized $K3$ surface of degree $2$.
Then the Stein factorization of $\Phi_{|H|}$ is of the form
\begin{equation*}
X\;\;\maprightsp{\rho}\;\; Y\;\; \maprightsp{\pi}\;\; \P^2,
\end{equation*}
where $\pi: Y\to \Pt$ is a finite double cover.
We have $h^0 (X, \OOO_X (mH))=m^2+2$ for every $m\ge 1$ 
by~\cite[Proposition~0.1]{Nikulin}.
Therefore the finite double cover $\pi: Y\to \Pt$  is defined by the equation
\begin{equation}\label{eq:CG}
w^2 + w \,C(x_0, x_1, x_2) + G(x_0, x_1, x_2)=0
\end{equation}
in the total space of the line bundle $V\to \Pt$ corresponding to 
the invertible sheaf $\OOO_{\Pt} (3)$,
where $w$ is a fiber coordinate of $V$,
$[x_0: x_1: x_2]$ is a homogeneous coordinate system of $\Pt$,
and $C$ and $G$ are homogeneous polynomials of degree $3$ and $6$ that are
regarded as sections of 
$V$ and $V\sp{\otimes 2}$, respectively.
If $C\ne 0$, then $\pi$ is separable,
while if $C=0$, then $\pi$ is purely inseparable.
\begin{definition}
An irreducible curve $F\st X$ is called a \emph{half-line of $(X, H)$}
if $FH=1$ holds.
A line $L\st\Pt$ is said to \emph{be splitting in $(X, H)$}
if the proper transform of $L$ in $X$ is non-reduced or reducible,
or equivalently,
if the scheme-theoretic pre-image $\pi\inv (L)\subset Y$ of $L$ by $\pi$
is non-reduced or reducible.
\end{definition}
Let $F$ be a half-line of $(X, H)$.
Then $\Phi_{|H|}$ induces an isomorphism from $F$ to a line $L\st \Pt$,
and this line $L$ is splitting in $(X, H)$.
In particular,
a half-line is a $(-2)$-curve. 
\begin{definition}
If $L\subset \Pt$ is a line splitting in $(X, H)$,
then the proper transform of $L$ in $X$ is written as $F+F\sprime$,
where $F$ and $F\sprime$ are half-lines of $(X, H)$.
These half-lines are said to be 
\emph{lying over $L$}.
We say that $L$ is   \emph{of non-reduced type}
if $F=F\sprime$,
while $L$ is \emph{of reduced type} if $F\ne F\sprime$.
\end{definition}
\begin{lemma}\label{lem:atmost3}
Suppose that $\pi$ is separable.
Then the number of splitting lines of non-reduced type is at most $3$.
\end{lemma}
\begin{proof}
Let $L$ be a splitting line of non-reduced type.
We choose  homogeneous coordinates $[x_0:x_1:x_2]$
of $\Pt$ such that $L$ is defined by $x_2=0$.
Putting $x_2=0$ in the defining equation~\eqref{eq:CG},
we see that the curve defined by 
\begin{equation}\label{eq:gamma}
w^2+w\,C(x_0, x_1, 0)+G(x_0, x_1, 0)=0
\end{equation}
in the total space of the line bundle $V|_L \to L$ on $L$ is non-reduced.
Let $\gamma(w, x_0, x_1)$ be the left-hand side of~\eqref{eq:gamma}.
Since $\ch k=2$,
we have $\partial \gamma/\partial w=C(x_0, x_1, 0)$.
Therefore $C(x_0, x_1, 0)$ is constantly equal to zero.
Thus we have shown  that the defining equation of a splitting line of non-reduced type
divides $C(x_0, x_1, x_2)$. Therefore, if $C\ne 0$,
then  the number of  splitting lines of non-reduced type is at most $\deg C=3$.
\end{proof}
Next we investigate the case where $\pi$ is purely inseparable.
In this case,  $\pi$ is given by the equation
\begin{equation}\label{eq:G}
w^2+G(x_0, x_1, x_2)=0.
\end{equation}
Note that
every splitting line is now of non-reduced type.
\begin{remark}\label{rem:Gamma}
Let $\Gamma(x_0, x_1, x_2)$ be a homogeneous polynomial of degree $3$.
Then the equations $w^2=G$ and $w^2=G+\Gamma^2$ define surfaces 
isomorphic over $\Pt$.
\end{remark}
We have the following relation between  splitting lines and  rational double points
of $Y$.
See~\cite{Artin1977}~or~\cite{GK} for the normal form of defining equations of rational double points
in characteristic $2$.
\begin{lemma}\label{lem:Q}
Let  $L\st \Pt$ be a line defined by $\ell(x_0, x_1, x_2)=0$.

{\rm (1)}
The line $L$ is splitting in $(X, H)$ if and only if 
there exist homogeneous polynomials $Q(x_0, x_1, x_2)$ and
$\Gamma (x_0, x_1, x_2)$ of degree $5$ and $3$, respectively,  such that
$G=\ell Q+\Gamma^2$.

{\rm (2)}
Suppose that $L$ is splitting in $(X, H)$,
and let $Q$ be a polynomial of degree $5$
such that $G+\ell Q$ is a square of a cubic polynomial.
 We denote by $T\st\Pt$ the quintic curve defined by $Q=0$.
Let $p$ be a point of $L$, and  $P$  the point of $Y$ such that $\pi (P)=p$.
Then 
$P$ is a smooth point of $Y$ if and only if  $p\notin T$,
and $P$ is an $A_1$-singular point of $Y$ if and only if
$T$ intersects $L$ transversely at $p$.
\end{lemma}
\begin{proof}
We can assume that $\ell=x_2$.
Since the curve defined by 
$w^2+G(x_0, x_1, 0)=0$
in  $V|_L$ is non-reduced,
we see that $G(x_0, x_1, 0)$ is the square of a polynomial of degree $3$.
Hence the assertion (1) follows.
Let $(x, y)$ be an affine coordinate system
of $\Pt$ with the origin $p$ such that $L$ is defined by $y=0$. 
We write~\eqref{eq:G} as
$w^2=g(x, y)$. Let $g_{ij}$ be the coefficient of $x^i y^j$ of $g$.
Then 
$P$ is a smooth point of $Y$ if and only if  $g_{01}\ne 0$ or $g_{10}\ne 0$,
and $P$ is an $A_1$-singular point of $Y$ if and only if
$g_{01}=g_{10}=0$ and $g_{11}\ne 0$.
Let $q(x, y)$ be the inhomogeneous polynomial corresponding to $Q$,
and let $q_{ij}$ be the coefficients of $x^i y^j$ of $q$.
Then, up to a multiplicative constant,
 we have $g_{01}=q_{00}$, $g_{10}=0$, $g_{11}=q_{10}$.
Therefore the assertion (2) follows.
\end{proof}
\begin{remark}
The polynomials $Q$ and $\Gamma$ such that
$G=\ell Q+\Gamma^2$ are not determined uniquely
by $G$ and $\ell$.
However, the homogeneous polynomial $Q|L$ on the line $L$ is
determined uniquely by $G$ and $\ell$.
\end{remark}
\section{Schr\"oer's Kummer surfaces as Zariski surfaces}\label{sec:schroeer}
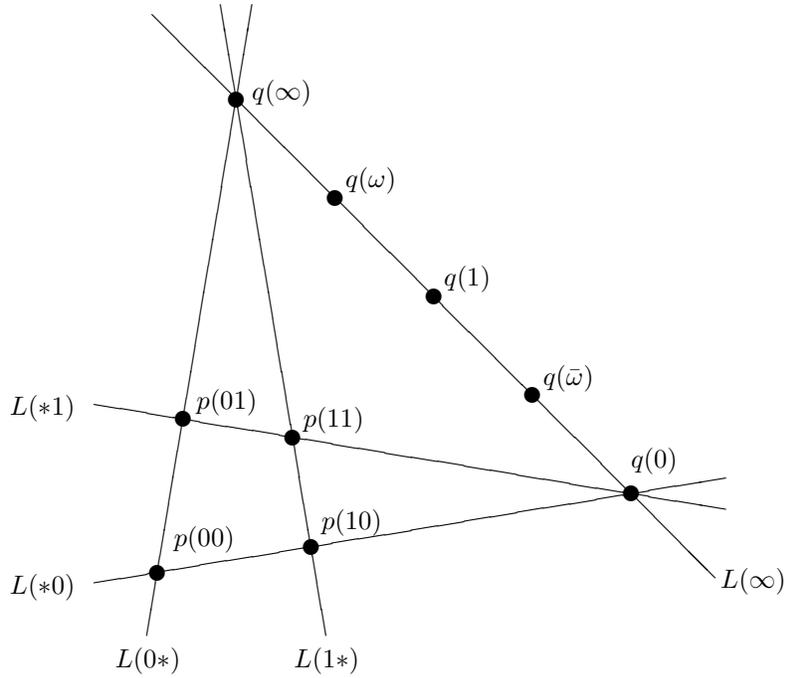
\begin{figure}\begin{center}
\setlength{\unitlength}{1.4mm}
\begin{picture}(70, 65)(0,0)
\put(10, 5){\line(1,6){10}}
\put(7, 2){$L(0*)$}
\put(27, 5){\line(-1, 6){10}}
\put(24, 2){$L(1*)$}
\put(5, 10){\line(6,1){60}}
\put(-3, 9){$L(*0)$}
\put(5, 27){\line(6, -1){60}}
\put(-3, 26){$L(*1)$}
\put(10.5, 64){\line(1, -1){53.5}}
\put(64.5,9.5){$L(\infty)$}
\put(18.5, 56){\circle*{1.4}}
\put(20, 56){$q(\infty)$}
\put(56, 18.5){\circle*{1.4}}
\put(56, 21){$q(0)$}
\put(37.25, 37.25){\circle*{1.4}}
\put(38.25, 38.25){$q(1)$}
\put(11,11){\circle*{1.4}}
\put(12.6,13.5){$p(00)$}
\put(23.8571, 23.8571){\circle*{1.4}}
\put(24.8571, 24.8571){$p(11)$}
\put(25.5946, 13.4324){\circle*{1.4}}
\put(26.5946, 15.1324){$p(10)$}
\put(13.4324, 25.5946){\circle*{1.4}}
\put(14.7324, 26.5946){$p(01)$}
\put(27.875,46.625){\circle*{1.4}}
\put(28.875,47.625){$q(\omega)$}
\put(46.625,27.875){\circle*{1.4}}
\put(47.625,28.875){$q(\bar\omega)$}
\end{picture}
\end{center}
\caption{The configuration of  splitting lines}\label{fig:config}
\end{figure}

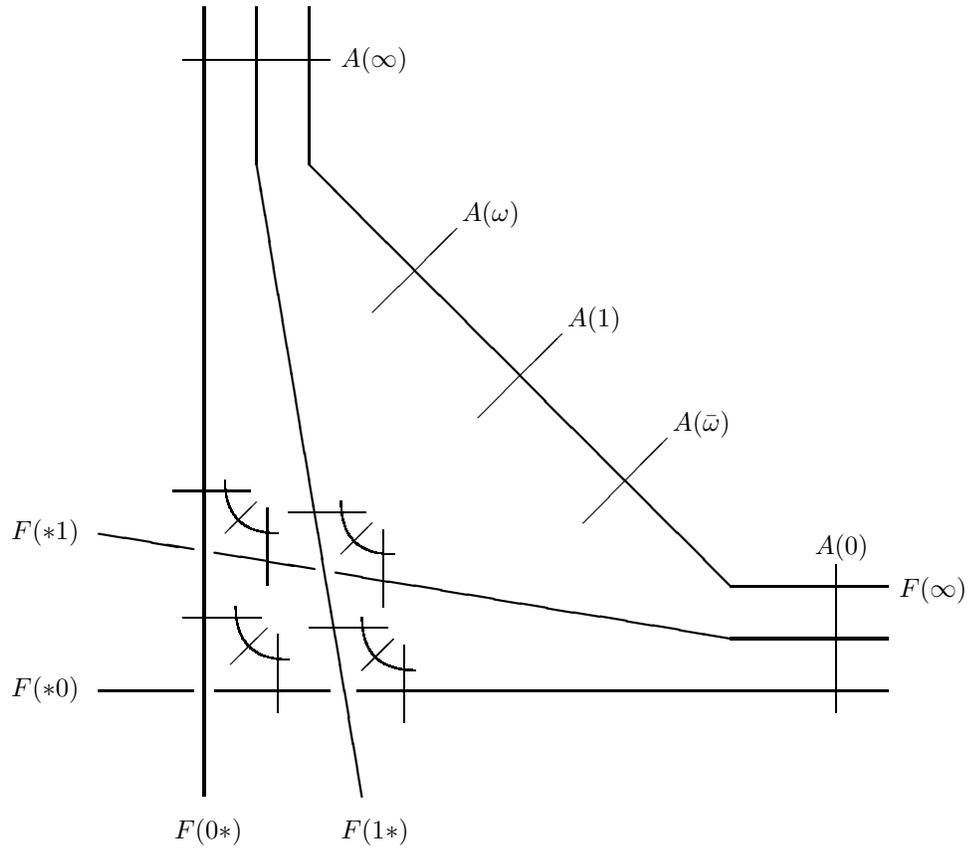
\begin{figure}
\begin{center}
\setlength{\unitlength}{1.4mm}
\begin{picture}(85, 85)(0,0)
\thinlines
\put(80, 18){\line(0,1){14}} \put(78, 33){$A(0)$}
\put(18, 80){\line(1,0){14}} \put(33, 79.3){$A(\infty)$}
\put(36,56){\line(1,1){8}} \put(44.5, 64.7){$A(\omega)$}
\put(46,46){\line(1,1){8}} \put(54.5, 54.7){$A(1)$}
\put(56,36){\line(1,1){8}} \put(64.5, 44.7){$A(\bar\omega)$}
\qbezier(23,28)(23,23)(28, 23)
\put(22.5, 22.5){\line(1,1){3.5}}
\put(18, 27){\line(1,0){7.4}}
\put(27, 18){\line(0,1){7.4}}
\qbezier(35,27)(35,22)(40, 22)
\put(35,22){\line(1,1){3}}
\put(35, 22){\line(1,1){3}}
\put(30, 26){\line(1,0){7.4}}
\put(39, 17){\line(0,1){7.4}}
\qbezier(27,35)(22,35)(22,40)
\put(22,35){\line(1,1){3}}
\put(17, 39){\line(1,0){7.4}}
\put(26, 30){\line(0,1){7.4}}
\qbezier(33,38)(33,33)(38, 33)
\put(33,33){\line(1,1){3}}
\put(28, 37){\line(1,0){7.4}}
\put(37, 28){\line(0,1){7.4}}
\thicklines
\put(10,20){\line(1,0){9}}
\put(21,20){\line(1,0){11}}
\put(34.5,20){\line(1,0){50.5}}
\put(1.7, 19.5){$F(*0)$}
\put(10,35){\line(6,-1){9}}
\put(21,33.167){\line(6,-1){9.5}}
\put(32.5,31.25){\line(6,-1){37.5}}
\put(70,25){\line(1,0){15}}
\put(1.7, 34.5){$F(*1)$}
\put(70,30){\line(1,0){15}}
\put(70,30){\line(-1,1){40}}
\put(30,70){\line(0,1){15}}
\put(86, 29){$F(\infty)$}
\put(20,10){\line(0,1){75}}
\put(17.2, 6){$F(0*)$}
\put(35,10){\line(-1,6){10}}
\put(25,70){\line(0,1){15}}
\put(33, 6){$F(1*)$}
\end{picture}
\end{center}
\caption{The configuration of half-lines and exceptional curves}\label{figure:curves}
\end{figure}

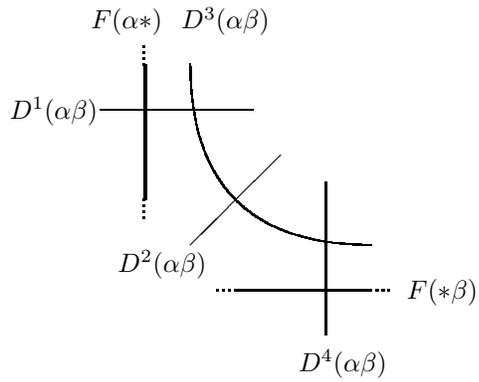
\begin{figure}
\begin{center}
\setlength{\unitlength}{1.2mm}
\begin{picture}(45, 45)(0,4)
\thinlines
\put(10,35){\line(1,0){17}}
\put(0, 34){$D^1(\alpha\beta)$}
\put(35,10){\line(0,1){17}}
\put(32,6){$D^4(\alpha\beta)$}
\put(30,30){\line(-1,-1){10}}
\put(12, 17){$D^2(\alpha\beta)$}
\qbezier(20, 40)(20,20)(40,20)
\put(19, 44){$D^3(\alpha\beta)$}
\thicklines
\qbezier[3](14.85, 40)(14.85, 41)(14.85, 42)
\put(15, 40){\line(0, -1){15}}
\qbezier[3](14.85, 25)(14.85, 24)(14.85, 23)
\put(9, 44){$F(\alpha*)$}
\qbezier[3](23, 15)(24, 15)(25, 15)
\put(40, 15){\line(-1,0){15}}
\qbezier[3](40, 15)(41, 15)(42, 15)
\put(44,14){$F(*\beta)$}
\end{picture}
\end{center}
\caption{The exceptional curves over $P(\alpha\beta)$}\label{figure:overD4}
\end{figure}

Let $r$ and $s$ be constants in $k$ such that $r\ne 0$, $s\ne 0$ and $r^3\ne s^3$.
Then Schr\"oer's supersingular $K3$ surface $X\rs$ 
defined in Introduction
is of Artin invariant $2$.
By Proposition~6.2~of~\cite{Sch}, 
the quotient surface $(C\times C)/\alpha_2$
of the $\alpha_2$-action on $C\times C$ defined by 
the vector field~\eqref{eq:delta} contains an open subset $U$ isomorphic to
$$
\Spec k[a,b,c]/(\,c^2+a(b^4+s^2 b^2)+b(a^4+r^2 a^2)\,).
$$
The singular locus of $U$ consists of four $D_4$-singular points
coming from the fixed points of the $\alpha_2$-action on the smooth part of $C\times C$.
Let 
$$
\pi\rs: Y\rs\to\Pt
$$
 be the purely inseparable double cover 
defined by
$$
w^2= [x_0 (x_1^4 + s^2 x_1^2 x_2^2) + x_1 (x_0^4 +r^2 x_0^2 x_2^2)] x_2,
$$
which is a projective completion of $U$.
Then  $Y\rs$  is birational to $X\rs$, 
and hence there exists a morphism 
$\rho\rs: X\rs\to Y\rs$ that is the minimal resolution.
The pull-back of a line of $\Pt$ by $\pi\rs\circ\rho\rs$ is a polarization $H\rs$
of degree $2$ of $X\rs$.
Then
$$
X\rs\;\;\maprightsp{\rho\rs}\;\; Y\rs\;\; \maprightsp{\pi\rs}\;\; \P^{2}
$$
is the Stein factorization of $\Phi_{|H\rs|}$.
The singular locus of  $Y\rs$ consists of four $D_4$-singular points
$P(00), P(01), P(10), P(11)$ in $U$ and five $A_1$-singular points
$Q(0), Q(1),  Q(\omega), Q({\bar\omega}), Q(\infty)$ lying on the line defined by $x_2=0$.
Here $\omega$ is a primitive third root of $1$, and $\bar\omega=\omega^2$.
These singular points are indexed in such a way that 
their images by $\pi\rs$ are given  in Table~\ref{table:points},
where $p(\ab):=\pi\rs (P(\ab))$ 
for $\ab=00,01,10,11$,
and $q(\gamma):=\pi\rs (Q(\gamma))$ 
for $\gamma= 0,1,\omega,\bar\omega,\infty$.
\begin{table}
$$
\begin{array}{ccccc}
p(00)&=&[0:0:1] \\
p(01)&=&[0:s:1] \\
p(10)&=&[r:0:1] \\
p(11)&=&[r:s:1] \\
\phantom{\pi\rs(P(00))} &&
\end{array}
\quad
\left|
\quad
\begin{array}{ccccc}
q(0)&=& [1:0:0]\\
q(1)&=& [1:1:0]\\
q(\omega)&=& [1:\omega:0]\\
q(\bar\omega)&=& [1:\bar\omega:0]\\
q(\infty)&=& [0:1:0]
\end{array}
\quad\right.
$$
\caption{The coordinates of the singular points of $Y\rs$}\label{table:points}
\end{table}
It is easy to see that the five lines listed below are splitting in $(X\rs, H\rs)$:
$$
\begin{array}{ll}
L(\infty):=\{x_2=0\}, & \\
L(0*):=\{x_0=0\}, &
L(1*):=\{x_0+r x_2=0\}, \\
L(*0):=\{x_1=0\}, &
L(*1):=\{x_1+s x_2=0\}. 
\end{array}
$$
To simplify the notation, we put
$$
\Pp:=\{00,01,10,11\},
\quad
\Qp:=\{0,1,\omega,\bar\omega,\infty\},
\quad
\SP:=\{\infty, 0*,1*,*0,*1\}.
$$
The configuration of the splitting lines $L(\lambda)$ ($\lambda\in\SP$)
and the points $p(\ab)$ ($\ab\in\Pp$), $q(\gamma)$ ($\gamma\in\Qp$) are given in Figure~\ref{fig:config}.
For a splitting line $L(\lambda)$ $(\lambda\in \SP)$,
we denote by $F(\lambda)$ the half-line of $(X\rs, H\rs)$ 
lying over $L(\lambda)$.
By blowing up $Y\rs$ at their singular points explicitly,
we see that the half-lines $F(\lambda)$ and the exceptional divisors of 
$\rho\rs:X\rs\to Y\rs$ intersect as in Figure~\ref{figure:curves}.
We denote the exceptional curves over the $D_4$-singular points
$P(\ab)$ $(\ab\in\Pp)$ as in Figure~\ref{figure:overD4},
and denote the exceptional curves over the $A_1$-singular points $Q(\gamma)$
$(\gamma\in \Qp)$ by $A(\gamma)$.
The polarized $K3$ surface $(X\rs, H\rs)$ has maximal rational double points.
Consider the sublattice
\begin{equation}\label{eq:latrs}
\lat\rs:=\lat_{\pol} \oplus 
\bigoplus\sb{\ab\in \Pp} \lat_{P(\ab)} \oplus 
\bigoplus\sb{\gamma\in \Qp}\lat_{Q(\gamma)}
\end{equation}
of $\NS (X\rs)$ with finite index,
as in the subsection~\ref{subsec:polmax}.
The lattice $\lat_{\pol}$ is of rank $1$ 
 generated by $h:=[H\rs]$,
and $\lat\dual_{\pol}$ is generated by $h\dual:=h/2$.
The lattice $\lat_{P(\ab)}$ is of rank $4$ with basis
$d^i(\ab):=[D^i(\ab)]$ ($i=1, \dots, 4$).
We denote the basis of $\lat_{P(\ab)}\dual$ dual to $d^1(\ab), \dots, d^4(\ab)$
by $d^1(\ab)\dual, \dots, d^4(\ab)\dual$.
The relation between $d^1(\ab), \dots, d^4(\ab)$ and $d^1(\ab)\dual, \dots, d^4(\ab)\dual$
are given by~\eqref{eq:dualD4}.
The lattice $\lat_{Q(\gamma)}$ is of rank $1$ generated by $a(\gamma):=[A(\gamma)]$,
and $\lat_{Q(\gamma)}\dual$ is generated by $a(\gamma)\dual:=-a(\gamma)/2$.
From Figures~\ref{figure:curves}~and~\ref{figure:overD4},
we see that the classes of half-lines $F(\lambda)$ ($\lambda\in\SP$)
are given as follows:
\begin{equation}\label{eq:HLclasses}
\begin{array}{ccl}
[F(\infty)]&=&
h\dual + a(0)\dual + a(\omega)\dual + a(1)\dual + a(\bar\omega)\dual + a(\infty)\dual, \cr
[F(0*)]&=&
h\dual + d^1 (00)\dual +d^1 (01)\dual + a(\infty)\dual, \cr
[F(1*)]&=&
h\dual + d^1 (10)\dual +d^1 (11)\dual + a(\infty)\dual, \cr
[F(*0)]&=&
h\dual + d^4 (00)\dual +d^4 (10)\dual + a(0)\dual, \cr
[F(*1)]&=&
h\dual + d^4 (01)\dual +d^4 (11)\dual + a(0)\dual.
\end{array}
\end{equation}
We then put
\begin{equation}\label{eq:dgrrs}
\dgr\rs:=(\lat\rs)\dual/\lat\rs =\dgr_{\pol} \oplus 
\bigoplus\sb{\ab\in \Pp} \dgr_{P(\ab)} \oplus 
\bigoplus\sb{\gamma\in \Qp}\dgr_{Q(\gamma)},
\end{equation}
which  is an $\F_2$-vector space of 
dimension $14$.
Since the discriminant of $\NS (X\rs)$ is $-2^{2\sigma (X\rs)}=-2^4$,
we see that $\NS(X\rs)/\lat\rs\st \dgr\rs$ is a subspace of dimension $5$.
It is easy to prove that the five elements
$$
\ol{[F(\lambda)]}:=[F(\lambda)] \bmod \lat\rs \quad(\lambda\in \SP)
$$
of $\NS(X\rs)/\lat\rs$ are linearly independent.
Therefore  $\NS (X\rs)$ is generated by the classes 
$h=[\pol\rs]$, $d^i(\ab)=[D^i(\ab)]$,
$a(\gamma)=[A(\gamma)]$ and $[F(\lambda)]$.
\begin{remark}
Suppose that $r^3=s^3$.
Then there exists $c\in \F_4\sptimes =\{1,\omega,\bar\omega\}$
such that $s=cr$ holds.
The three points $p(00), p(11)$ and $ q(c)$
on $\Pt$ are collinear.
Let $M$ be the line passing through these points.
Then $M$ is a splitting line for $(X_{r, cr}, H_{r, cr})$.
Let $G$ be the half-line lying over $M$.
By blowing up $Y_{r, cr}$ at the points $P(00), P(11)$ and $ Q(c)$,
we see that
$$
[G]=h\dual +d^2(00)\dual+ d^2(11)\dual +a(c)\dual.
$$
Note that $d^2(\ab)\dual \equiv d^1(\ab)\dual +d^4(\ab)\dual \bmod \lat_{P(\ab)}$ 
by~\eqref{eq:dualD4}.
Hence $\ol{[G]}:=[G]\bmod \lat_{r, cr}$ is linearly independent from 
the set of vectors 
$\ol{[F(\lambda)]}\; (\lambda\in \LLL)$
in $\dgr_{r,cr}$.
In particular, the linear subspace $\NS(X_{r, cr})/\lat_{r, cr}$ 
of $\dgr_{r,cr}$ is of dimension $6$
generated by $\ol{[F(\lambda)]}\; (\lambda\in \LLL)$ and $\ol{[G]}$, 
and the Artin invariant of $X_{r, cr}$ is $1$.
\end{remark}
\section{Proof of main theorem}\label{sec:proof}
Note that supersingular $K3$ surfaces with Artin invariant $1$
are isomorphic to each other~\cite{DK, Ogus}.
Therefore
it is enough to prove Theorem~\ref{thm:main} 
under the additional assumption that the Artin invariant of $X\sprime$ is $2$.
\par
\medskip
We choose a Schr\"oer's Kummer surface $X$ with $\sigma (X)=2$.
To fix the ideas, we choose $s\in k\sm \F_4$, and  put $X:=X_{1,s}$, and set 
$$
H:=H_{1,s}, \quad
Y:=Y_{1,s},\quad
\rho:=\rho_{1,s},\quad
\pi:=\pi_{1,s}, \quad
\dgr:=\dgr_{1,s}, \quad
\lat:=\dgr_{1,s}.
$$
Let $X\sprime$ be a supersingular $K3$ surface with Artin invariant $2$.
Theorem~\ref{thm:RS} implies that
$\NS(X)$ and $\NS(X\sprime)$ are isomorphic.
By~Proposition~\ref{prop:nef}, there exists an isomorphism $\phi:\NS (X)\isom \NS (X\sprime)$
such that $\phi\otimes\R$ maps $\Nef (X)$ to $\Nef (X\sprime)$.
We fix such an isomorphism $\phi$ once and for all.
By~Corollary~\ref{cor:pol},
we have a polarization $H\sprime$ of $X\sprime$ with degree $2$ such that
$[H\sprime]=\phi ([H])$.
As before,
let 
\begin{equation*}
X\sprime\;\;\maprightsp{\rho\sprime}\;\; Y\sprime\;\; \maprightsp{\pi\sprime}\;\; \P^{2}
\end{equation*}
be the Stein factorization of $\Phi_{|H\sprime|}$.
By~Corollary~\ref{cor:commutative},
there exist bijections
$\phi_{\EEE}: \EEE\isom \EEE\sprime$ and 
$\phi_{\Sing}: \Sing (Y)\isom \Sing(Y\sprime)$
such that the  diagram~\eqref{eq:commutative} is commutative.
For  $P\in \Sing (Y)$, we write $P\sprime\in \Sing (Y\sprime)$
instead of $\phi_{\Sing} (P)$,
and for  $E\in \EEE$,
we write $E\sprime\in \EEE\sprime$ instead of  $\phi_{\EEE} (E)$.
Therefore $\Sing (Y\sprime)$ consists of four $D_4$-singular points 
$P(\ab)\sprime$ ($\ab\in \Pp$), 
and five $A_1$-singular points 
$Q(\gamma)\sprime$ ($\gamma\in \Qp$).
For example,
the $(-2)$-curves contracted to $P(\ab)\sprime$ by $\rho\sprime$
are $D^1(\ab)\sprime$, $D^2(\ab)\sprime$, $D^3(\ab)\sprime$ and $D^4(\ab)\sprime$.
We then put 
$$
p(\ab)\sprime:=\pi\sprime (P(\ab)\sprime)
\quand
q(\gamma)\sprime:=\pi\sprime (Q(\gamma)\sprime).
$$
We also set
\begin{equation}\label{eq:latsprime}
\lat\sprime:=
\lat_{\pol\sprime} \oplus 
\bigoplus\sb{\ab\in \Pp} \lat_{P(\ab)\sprime} \oplus 
\bigoplus\sb{\gamma\in \Qp}\lat_{Q(\gamma)\sprime}
\end{equation}
and
\begin{equation}\label{eq:dgrsprime}
\dgr\sprime:=(\lat\sprime)\dual/\lat\sprime =
\dgr_{\pol\sprime}\oplus
\bigoplus\sb{\ab\in \Pp} \dgr_{P(\ab)\sprime} \oplus 
\bigoplus\sb{\gamma\in \Qp}\dgr_{Q(\gamma)\sprime}
\end{equation}
as~\eqref{eq:latrs}~and~\eqref{eq:dgrrs}. 
Note that $\phi$ induces isomorphisms
$$
\phi_\lat:\lat\isom\lat\sprime\quand
\phi_\dgr:\dgr\isom\dgr\sprime
$$
that are compatible with the direct-sum 
decompositions~\eqref{eq:latrs}, \eqref{eq:latsprime}, 
and \eqref{eq:dgrrs}, \eqref{eq:dgrsprime}.
\par
\medskip
Let $L\st\Pt$ be a line splitting in $(X, H)$,
and let $F$ be the half-line of $(X, H)$ lying over  $L$.
We can define  a line $L\sprime\st\Pt$
splitting in $(X\sprime, H\sprime)$
and a half-line $F\sprime$ of $(X\sprime, H\sprime)$ lying over  $L\sprime$ as follows.
\begin{claim}\label{claim:uniqueDsprime}
There exists a unique effective divisor $D\sprime$
that represents $\phi([F])$.
\end{claim}
\begin{proof}
Since $\phi([F])^2=-2$ and $\phi([F])[H\sprime]=1$,
there exists an effective divisor $D\sprime$
that represents $\phi([F])$.
Let $D\sprime=\Gamma\sprime+M\sprime$ be the decomposition
of $D\sprime$ into the sum of the fixed part $\Gamma\sprime$ and the movable part $M\sprime$.
Suppose that $M\sprime\ne 0$.
If  $M\sprime H\sprime=0$,
then $M\sp{\prime \;2}<0$ because $[H\sprime]\sperp$ is negative-definite.
Therefore we have  $M\sprime H\sprime>0$.
Since $\Gamma\sprime H\sprime\ge 0$,
we have $M\sprime H\sprime=1$,
which implies that $\Phi_{|H\sprime|}$ induces an isomorphism
from $M\sprime$ to a line on $\Pt$.
Hence $M\sprime$ is a smooth rational curve, which is a contradiction.
\end{proof}
Since $D\sprime H\sprime=1$,
there exists a unique irreducible component $F\sprime$ of $D\sprime$
such that $F\sprime H\sprime=1$.
Then $F\sprime$ is a half-line of $(X\sprime, H\sprime)$.
We define  $L\sprime\st\Pt$ to be  the image of $F\sprime$ by $\rho\sprime\circ \pi\sprime$.
\begin{claim}\label{claim:codedeterminesL}
Let $F\spprime$ be a half-line for $(X\sprime, H\sprime)$
lying over $L\sprime$.
Then $\ol{[F\spprime]}=\ol{[F\sprime]}$ holds in $\dgr\sprime$,
where
$\ol{[F\spprime]}=[F\spprime]\;\bmod\lat\sprime$
and $\ol{[F\sprime]}=[F\sprime]\;\bmod\lat\sprime$.
\end{claim}
\begin{proof}
The case where $F\sprime=F\spprime$ is obvious.
Suppose that $F\sprime\ne F\spprime$.
Then $F\sprime+F\spprime$ is the total transform of $L\sprime$ in $X\sprime$
minus a linear combination of curves in $\EEE\sprime$,
and hence $[F\sprime]+[F\spprime]\in \lat\sprime$.
Because $\dgr\sprime$ is a $2$-elementary abelian group, we obtain
$\ol{[F\spprime]}=\ol{[F\sprime]}$.
\end{proof}
\begin{claim}\label{claim:code} 
We have $\phi_{\dgr} (\ol{[F]})=\ol{[F\sprime]}$.
\end{claim}
\begin{proof}
Since $\phi([F])=[D\sprime]$,
we have $\phi_{\dgr} (\ol{[F]})=\ol{[D\sprime]}$.
Since $D\sprime-F\sprime$ is effective and $(D\sprime-F\sprime)H\sprime=0$,
each irreducible component of $D\sprime-F\sprime$ is contracted to a point by 
$\rho\sprime$.
Therefore we have $[D\sprime]-[F\sprime]\in \lat\sprime$,
and hence $\ol{[D\sprime]}=\ol{[F\sprime]}$.
\end{proof}
Now we have half-lines $F(\lambda)\sprime$ and splitting lines $L(\lambda)\sprime$
of $(X\sprime, H\sprime)$ for each $\lambda\in \LLL$.
By Claim~\ref{claim:code},
the elements $\ol{[F(\lambda)\sprime]}$ of $\dgr\sprime$ are distinct to each other.
Hence, by Claim~\ref{claim:codedeterminesL},
the lines $L(\lambda)\sprime$ are distinct to each other.
\begin{claim}\label{claim:passingthrough}
Let $P$ be a point of $\Sing (Y)$.
If $\pi (P) \in L(\lambda)$, then $\pi\sprime(P\sprime)\in L(\lambda)\sprime$.
\end{claim}
\begin{proof}
If $\pi (P) \in L(\lambda)$, 
then the $P$-component of $\ol{[F(\lambda)]}\in \dgr$ is not zero by~\eqref{eq:HLclasses}.
Hence the $P\sprime$-component of $\ol{[F(\lambda)\sprime]}\in \dgr\sprime$ is not zero 
by Claim~\ref{claim:code}.
Consequently, there exists $E\sprime\in \EEE\sprime_{P\sprime}$ such that
$F(\lambda)\sprime E\sprime\ne 0$.
Therefore the image $L(\lambda)\sprime$ of $F(\lambda)\sprime$ passes through $\pi\sprime(P\sprime)\in
\Pt$.
\end{proof}
\begin{claim}\label{claim:nonreduced}
The splitting line $L(\lambda)\sprime$ is of non-reduced type
for any $\lambda\in \SP$.
\end{claim}
\begin{proof}
Let $G\sprime$ be an arbitrary half-line of $(X\sprime, H\sprime)$
lying over $L(\lambda)\sprime$.
Then the class $g\sprime:=[G\sprime]\in \NS(X\sprime)$ satisfies the following:
\begin{itemize}
\item[(i)] $(g\sprime)^ 2=-2$,
\item[(ii)] $g\sprime [H\sprime]=1$, and 
\item[(iii)] for every $E\sprime\in \EEE\sprime$, we have $g\sprime [E\sprime]\ge 0$.
\end{itemize}
Suppose that $L(\lambda)\sprime$ is of reduced type.
Then there exists a half-line $F\spprime$lying over $L(\lambda)\sprime$
that is distinct  from
$F(\lambda)\sprime$.
Since $[F\spprime][F(\lambda)\sprime]\ge 0$,
we have $[F\spprime]\ne [F(\lambda)\sprime]$.
By Claim~\ref{claim:codedeterminesL}, we have 
$\ol{[F(\lambda)\sprime]}=\ol{[F\spprime]}$ in $\dgr\sprime$.
Consequently, it is enough to show that there exists only one class $g\sprime$
in $\NS (X\sprime)$ satisfying (i),  (ii), (iii)
 above and 
\begin{itemize}
\item[(iv)] $\ol{(g\sprime)}=\ol{[F(\lambda)\sprime]}=\phi_\dgr (\ol{[F(\lambda)]})$,
\end{itemize}
where the second equality follows from Claim~\ref{claim:code}.
We denote by $g\sprime_{H\sprime}$ and $g\sprime_{P\sprime}$ the
$H\sprime$- and $P\sprime$-components of $g\sprime$, respectively,
where $P\sprime\in \Sing (S\sprime)$.
By (ii),
we have $g\sprime_{H\sprime}=[H\sprime]/2$. 
Combining this with (i), we have
\begin{equation}\label{eq:5/2}
\sum_{\ab\in \Pp} (g\sprime_{P(\ab)\sprime})^2 +
\sum_{\gamma\in\Qp} (g\sprime_{Q(\gamma)\sprime})^2
=-5/2.
\end{equation}
\par
{\it The case where $\lambda=\infty$.}
By (iii), (iv), \eqref{eq:HLclasses} and Lemmas~\ref{lem:A},~\ref{lem:D},
we have
$$
(g\sprime_{P(\ab)\sprime})^2=0 \;\;\textrm{or}\;\; \le -2
\quand 
 (g\sprime_{Q(\gamma)\sprime})^2=-1/2 \;\;\textrm{or}\;\; \le -9/2.
$$
Combining this with~\eqref{eq:5/2}, we have
$$
(g\sprime_{P(\ab)\sprime})^2=0
\quand 
 (g\sprime_{Q(\gamma)\sprime})^2=-1/2.
$$
By (iii) and Lemmas~\ref{lem:A},~\ref{lem:D} again, we have 
$$
g\sprime_{P(\ab)\sprime}=0
\quand 
 g\sprime_{Q(\gamma)\sprime}=-[A(\gamma)\sprime]/2.
$$
Thus the uniqueness of $g\sprime$ is proved.
\par
{\it The case where $\lambda=0*$.}
By (iii), (iv), \eqref{eq:HLclasses} and Lemmas~\ref{lem:A},~\ref{lem:D},
we have
$$
\begin{array}{rll}
(g\sprime_{P(\ab)\sprime})^2&=-1 \;\;\textrm{or}\;\; \le -3& \;\;\textrm{if $\ab=00$ or $01$,}\\
(g\sprime_{P(\ab)\sprime})^2&=0 \;\;\textrm{or}\;\; \le -2& \;\;\textrm{if $\ab=10$ or $11$,}\\
(g\sprime_{Q(\gamma)\sprime})^2&=-1/2 \;\;\textrm{or}\;\; \le -9/2& \;\;\textrm{if $\gamma=\infty$,}\\
(g\sprime_{Q(\gamma)\sprime})^2&=0 \;\;\textrm{or}\;\; \le -2& \;\;\textrm{if $\gamma\ne\infty$.}
\end{array}
$$
Combining this with~\eqref{eq:5/2}, we have
$$
\begin{array}{l}
(g\sprime_{P(00)\sprime})^2=(g\sprime_{P(01)\sprime})^2=-1,\quad
(g\sprime_{P(10)\sprime})^2=(g\sprime_{P(11)\sprime})^2=0,\\
(g\sprime_{Q(\infty)\sprime})^2=-1/2,\quad
(g\sprime_{Q(\gamma)\sprime})^2=0 \quad\textrm{for $\gamma\ne \infty$}.
\end{array}
$$
By (iii) and Lemmas~\ref{lem:A},~\ref{lem:D} again, we have 
$$
\begin{array}{l}
g\sprime_{P(00)\sprime}=\delta^1(00), \quad g\sprime_{P(01)\sprime}=\delta^1(01),\quad
g\sprime_{P(10)\sprime}=g\sprime_{P(11)\sprime}=0,\\
g\sprime_{Q(\infty)\sprime}=-[A(\infty)\sprime]/2,\quad
g\sprime_{Q(\gamma)\sprime}=0 \quad\textrm{for $\gamma\ne \infty$},
\end{array}
$$
where 
$$
\delta^1(\ab)=-[D^1(\ab)\sprime] -[D^2(\ab)\sprime]/2-[D^3(\ab)\sprime]-[D^4(\ab)\sprime]/2.
$$
(See~\eqref{eq:dualD4}.)
Thus the uniqueness of $g\sprime$ is proved.
\par
The other cases $\lambda=1*,*0,*1$ can be treated in the same way.
\end{proof}
We have 
five distinct splitting lines
$L(\lambda)\sprime$ ($\lambda\in \SP$) for  $(X\sprime, H\sprime)$,
which are of non-reduced type by Claim~\ref{claim:nonreduced}.
By Lemma~\ref{lem:atmost3},
we see that $\pi\sprime:Y\sprime\to \Pt$ is purely inseparable.
By Claim~\ref{claim:passingthrough},
the configuration of the lines $L(\lambda)\sprime$ 
and the points $p(\ab)\sprime, q(\gamma)\sprime$
are exactly the same as the   configuration depicted in Figure~\ref{fig:config}
with superscript prime (${}\sprime$) being put to everything.
\par
\medskip
There exists a homogeneous coordinate system
$[x: y: z]$ of $\Pt$ such that
$$
\begin{array}{l}
q(\infty)\sprime=[0:1:0],\quad
q(1)\sprime=[1:1:0],\quad
q(0)\sprime=[1:0:0],\\
p(00)\sprime=[0:0:1],\quad
p(10)\sprime=[1:0:1].
\end{array}
$$
We put
$$
p(01)\sprime=[0:t:1],
$$
where $t$ is a non-zero constant.
Then we have $p(11)\sprime=[1:t:1]$ by Figure${}\sprime$~\ref{fig:config}.
Let
$$
w^2=G(x, y, z)
$$
be the defining equation of $Y\sprime$,
where $G$ is a homogeneous polynomial of degree $6$, 
and let $G_{lmn}$ ($l+m+n=6$) be the coefficient of $x^l y^m z^n$ in $G$.
By~Remark~\ref{rem:Gamma}, we can assume
\begin{equation*}\label{eq:lmneven}
G_{lmn}=0\;\;\;\;\textrm{if \;\;\;\;$l\equiv m\equiv n\equiv 0\;\bmod 2$}.
\end{equation*}
Because $L(0*)\sprime=\{x=0\}$ is splitting, Lemma~\ref{lem:Q}(1) implies
\begin{equation*}\label{eq:0star}
G_{015}=G_{033}=G_{051}=0.
\end{equation*}
Because $L(*0)\sprime=\{y=0\}$ is splitting, Lemma~\ref{lem:Q}(1) implies
\begin{equation*}\label{eq:star0}
G_{105}=G_{303}=G_{501}=0.
\end{equation*}
Because $L(\infty)\sprime=\{z=0\}$ is splitting, Lemma~\ref{lem:Q}(1) implies
\begin{equation*}\label{eq:infty}
G_{150}=G_{330}=G_{510}=0.
\end{equation*}
Therefore we have
$$
G(x, y, z)=x y z \,C(x, y, z),
$$
where $C$ is a homogeneous polynomial of degree $3$.
By Lemma~\ref{lem:Q}(2),
the line $L(0*)\sprime=\{x=0\}$ and the quintic curve defined by 
$y z \,C(x, y, z)=0$ intersect transversely at $q(\infty)\sprime$
and with multiplicity $\ge 2$ at $p(00)\sprime$ and $p(01)\sprime$.
Therefore there exists a constant $A$ such that
$y z C(0, y, z)=A y^2 z (y+tz)^2$.
In particular, we obtain
\begin{equation*}\label{eq:0star2}
G_{132}=G_{114}=0\quand G_{123}=t^2 G_{141}.
\end{equation*}
By Lemma~\ref{lem:Q}(2),
the line  $L(*0)\sprime=\{y=0\}$ and the  curve 
$x z \,C(x, y, z)=0$ intersect transversely at $q(0)\sprime$
and with multiplicity $\ge 2$ at $p(00)\sprime$ and $p(10)\sprime$.
Therefore there exists a constant $B$ such that
$x z C(x, 0, z)=B x^2 z (x+z)^2$.
In particular, we obtain
\begin{equation*}\label{eq:star02}
G_{312}=G_{114}=0\quand G_{213}= G_{411}.
\end{equation*}
By Lemma~\ref{lem:Q}(2),
the line  $L(\infty)\sprime=\{z=0\}$ and the  curve 
$x y \,C(x, y, z)=0$ intersect transversely at the five points $q(\gamma)\sprime$
($\gamma\in \Qp)$.
In particular, the curve $x y \,C(x, y, z)=0$ passes through $q(1)\sprime$,
and hence we obtain
\begin{equation*}\label{eq:inf2}
G_{141}+G_{231}+G_{321}+G_{411}=0.
\end{equation*}
Combining these, we see that $Y\sprime$ is defined  by
$$
w^2=x y z (t^2 a y z^2 + d x z^2 + a y^3 +b x y^2 + c x^2 y + d x^3), 
$$
where $a,b,c, d$ are constants such that $a+b+c+d=0$.
Because $L(1*)\sprime=\{x=z\}$ is splitting,
the polynomial $yz^2(t^2yz^2+ay^3+bzy^2+cz^2y)$ of $y$ and $z$  
is a square of a cubic polynomial.
Therefore $b=0$.
Because $L(*1)\sprime=\{y=tz\}$ is splitting,
the polynomial $tx z^2(dxz^2+bt^2xz^2+ctx^2z+dx^3)$ of $x$ and $z$  
is a square of a cubic polynomial.
Therefore $c=0$.
Because $a+b+c+d=0$, we have $a=d$.
Therefore $Y\sprime$ is defined by
$$
w^2=x y z (t^2  y z^2 +  x z^2 +  y^3 +  x^3).
$$
Hence $Y\sprime$ is isomorphic to Schr\"oer's  normal $K3$ surface $Y_{t, 1}$,
and hence $X\sprime$ is isomorphic to  Schr\"oer's  Kummer surface  $X_{t, 1}$.
\begin{remark}
In~\cite{phoshimada},
it is  shown that every  supersingular $K3$ surface in characteristic $5$
with Artin invariant $\le 3$ is obtained as a double cover of
the projective plane
with the branch curve defined by   $y^5-f(x)=0$,
where $f(x)$ is a polynomial of degree $6$,
and hence it is unirational.
\end{remark}
\bibliographystyle{plain}

\end{document}